\documentclass[twoside]{article}%
\usepackage{graphicx}
\usepackage{amsmath}
\usepackage{theorem}
\usepackage{amsfonts}
\usepackage{amssymb}%
\theoremstyle{change}
\newtheorem{theorem}{Theorem}[section]

\newtheorem{corollary}[theorem]{Corollary}

\newtheorem{lemma}[theorem]{Lemma}
\newtheorem{proposition}[theorem]{Proposition}
{
\theorembodyfont{\upshape}

\newtheorem{definition}[theorem]{Definition}
\newtheorem{example}[theorem]{Example}

\newtheorem{notation}[theorem]{Notation}

\newtheorem{remark}[theorem]{Remark}

}
\newenvironment{proof}[1][Proof]{\noindent \textbf{#1.} }{\ \rule{0.5em}{0.5em}}
\begin{document}

\title{Hamilton's injectivity radius estimate for sequences with almost nonnegative
curvature operators\thanks{To appear in \emph{Communications in Analysis and
Geometry.}}}
\author{Bennett Chow\thanks{Partially supported by NSF grant DMS-9971891.}\\University of California, San Diego
\and Dan Knopf\thanks{Partially supported by NSF grant DMS-0202796.}\\University of Wisconsin (current address:\ University of Iowa)
\and Peng Lu\\McMaster University (current address: University of Oregon)}
\date{September 17, 2001 (revised May 1, 2002)}
\maketitle

\section{ Introduction}

\footnotetext[1]{We would like to thank the referee for suggesting changes to
improve the presentation of the paper.}In recent years, one of the focuses in
the study of the Ricci Flow on Riemannian manifolds has been on classifying
the singularities that form in low dimensions. In particular, Hamilton has
obtained partial classification theorems in dimensions three and four. These
classifications are in the sense of obtaining and classifying pointed limits,
provided they exist, of dilations of a solution to the Ricci Flow about
sequences of points and times tending to the singularity time, after passing
to a suitable subsequence. In dimension four, Hamilton's classification
together with his geometric--topological surgery methods yield a
classification of diffeomorphism types of those compact $4$-manifolds with
positive isotropic curvature that do not admit any incompressible
$3$-dimensional space form not diffeomorphic to either $\mathbb{S}^{3}$ or
$\mathbb{R}\mathbb{P}^{3}$ \cite{H2}. (Note that \cite{MM} obtained an earlier
classification of homeomorphism types of compact simply-connected
$n$-manifolds with positive isotropic curvature by using harmonic maps.) In
dimension three, Hamilton's classification of singularities in \cite{H1} plays
a major role in his program for approaching Thurston's Geometrization
Conjecture \cite{T} by Ricci Flow methods. It is conjectured by Hamilton that
for the volume normalized Ricci flow on a compact $3$-manifold, after a finite
number of geometric--topological surgeries at some finite sequence of times,
the solution will exist for all time and the curvature will remain uniformly
bounded. If this is so, then the $3$-manifold admits a geometric decomposition
by \cite{H3}.

A fundamental tool used to obtain limits of sequences of solutions to the
Ricci Flow is the Gromov-type compactness theorem of Hamilton \cite{H4}. (See
\cite{P} for a survey of compactness theorems in Riemannian geometry.) As is
usual for compactness theorems, the main assumptions are bounded curvature
(which by the Bernstein--Shi estimates \cite{S} and \S 7 of \cite{H1} implies
bounds on all derivatives of curvature for solutions of the Ricci flow) and a
lower bound for the injectivity radius at a point (which by \cite{CLY} or
\cite{CGT} implies injectivity radii estimates at all points, depending on
distance). For the sequences arising from dilations about a singularity, the
curvature bound follows from the conditions imposed on the choice of points
and times and choice of dilation factor. (See \S 16 of \cite{H1}.) This leaves
one with the problem of obtaining an injectivity radius bound. In fact,
Hamilton's Little Loop Lemma (\S 15 of \cite{H1}) asserts that for a solution
of the Ricci Flow, a strengthened injectivity radius estimate should follow
from a suitable differential Harnack inequality of Li--Yau--Hamilton type. By
\cite{H5}, a differential Harnack inequality holds for complete solutions of
the Ricci flow with nonnegative curvature operator. See \cite{LY} for the
seminal result of this type for solutions of the heat equation. For Type I
singularities in dimension three, Hamilton proved an isoperimetric estimate
that implies an injectivity radius estimate (\S 23 of \cite{H1}). For Type II
singularities in dimension three, Hamilton also conjectures that there is an
injectivity radius estimate. In fact, the Little Loop Lemma, which is
conjectured to be true for all solutions of the Ricci Flow on compact
$3$-manifolds, subsumes this conjecture. Similarly, for K\"{a}hler manifolds
with positive bisectional curvature, there is an injectivity radius estimate
useful for the study of the K\"{a}hler--Ricci flow (see \cite{CT}). 

Although the above conjecture is still open, there is an important case where
an injectivity radius estimate should be true. Namely, consider a sequence of
complete (in practice, usually compact) solutions to the Ricci flow with
bounded curvature on a common time interval such that:

\begin{itemize}
\item The diameters are tending to infinity. Note that in the bounded diameter
case, either there is an injectivity radius estimate for the sequence or the
sequence collapses. In the former case, one can obtain a limit; and in the
latter case, when the sequence of underlying topological manifolds is a fixed
compact $3$-dimensional manifold $\mathcal{M}^{3}$, Cheeger--Gromov theory
proves that $\mathcal{M}^{3}$ is classified as a graph manifold.

\item The curvature operators are tending to nonnegative. (In dimension three,
this follows by \S 24 of \cite{H1} or \S 4 of \cite{H3} for a sequence of
solutions arising from dilations about a singularity.)

\item The curvatures at the origin are uniformly bounded from below by a
positive constant. (Otherwise we have the split case, which we hope to address
in a future article.)
\end{itemize}

\noindent In this case, Hamilton has claimed an injectivity radius estimate
after passing to a suitable subsequence. (See Theorem 25.1 in \cite{H1} and
Theorem \ref{InjEstForNNCO}, below.) One of the applications of Theorem 25.1
of \cite{H1} is in the proof of the classification of $4$-manifolds with
positive isotropic curvature given in \cite{H2}. (See all three subsections of
the \textquotedblleft Recovering the manifold from surgery\textquotedblright%
\ section of that paper.)

The purpose of this paper is to give a complete proof of the aforementioned
injectivity radius estimate. The reason for giving this new proof is related
to the possibility of collapse. In particular, there appears to be a gap in
the argument in \cite{H1}. We shall explain this in more detail in later
sections. (See the remarks before and after Example \ref{Collapse} and also
Remark \ref{TiltingFrames}.) The main overall structure of our proof is the
same as Hamilton's. However, our approach makes essential changes in the
construction of Busemann-type sublevel sets, and relies on new arguments to
establish the crucial fact that they are ultimately bounded. In particular,
our proof does not rely on the continuity of the function $\ell_{\infty}$
introduced in Lemma 25.3 of \cite{H1}. The main technical innovations of our
method are in sections \ref{mimicking} and \ref{boundedness}. We summarize
them here for the convenience of the reader:

\begin{itemize}
\item In \S 25 of \cite{H1}, it is argued that the distances to the cut loci
at the origins along any sequence of solutions to the Ricci flow satisfying
Definition \ref{almost nonnegative curvature operators} converge to a
continuous function $\ell_{\infty}:\mathcal{S}_{1}^{n-1}\rightarrow\left[
0,\infty\right]  $, where $\mathcal{S}_{1}^{n-1}$ is the standard $\left(
n-1\right)  $-sphere of radius $1$. There are subtle difficulties with this
approach relating to the possibility of collapse. (See Example \ref{Collapse}
and the Remark that follows.) To get around these difficulties, we define an
alternate function $\sigma_{\infty}:\mathcal{S}_{1}^{n-1}\rightarrow\left[
0,\infty\right]  $ in (\ref{sigma-infinity}) by a $\limsup$, which obviates
proving continuity.

\item In \S 25 of \cite{H1}, a set $\mathcal{D}$ of \emph{distinguished
directions} is defined as $\ell_{\infty}^{-1}\left(  \infty\right)  $, ``those
[directions] in which we can go off to infinity without hitting the cut
locus.'' We replace this by a set $\mathcal{R}_{\infty}\doteqdot\sigma
_{\infty}^{-1}\left(  \infty\right)  $ of \emph{ray-like directions.} The set
$\mathcal{R}_{\infty}$ is nonempty, and there exist arbitrarily long
minimizing geodesics in directions arbitrarily close to each of its members.
(See Remark \ref{TiltingFrames}.)

\item For robustness under the action of passing to subsequences, we find in
$\mathcal{R}_{\infty}$ for any $\varepsilon>0$ a finite subset $\left\{
V_{\alpha}\right\}  $ such that no member of $\mathcal{R}_{\infty}$ lies more
than distance $\varepsilon$ away from some $V_{\alpha}$, and (most
importantly) such that$\ $the $\limsup$ in Definition (\ref{sigma-infinity})
is attained as a limit for each $V_{\alpha}$.

\item In \S 25 of \cite{H1}, sets $N_{i}$ are defined that act as substitutes
for the sublevel sets of a Busemann function, and properties of the function
$\ell_{\infty}$ are invoked to claim that they are uniformly bounded. We
replace the $N_{i}$ by sets $N_{i}\left(  L\right)  \equiv N_{i}\left(
L,1\right)  $ that depend on a length scale $L\gg1$. (See Definition
(\ref{Busemann}) in \S \ref{mimicking}.) We do not show that the sets
$N_{i}\left(  L\right)  $ are uniformly bounded. But the key innovation in our
method is the proof of the boundedness property (Proposition
\ref{N_iLK-are-bounded}), which states, roughly, that by going far enough out
in the sequence, depending on $L$, one can bound the size of all remaining
$N_{i}\left(  L\right)  $ independently of $L$.
\end{itemize}

\paragraph{Acknowledgment}

This paper was prepared while the authors enjoyed the hospitality provided
during the summer of 2001 by the National Center for Theoretical Sciences in
Hsinchu, Taiwan. We wish to thank the NCTS for providing partial support and a
wonderful research environment.

\section{Hamilton's injectivity radius estimate}

We now recall the setup from \S 25 of \cite{H1}. Consider a sequence
\[
\left\{  \mathcal{M}_{i}^{n},g_{i}\left(  t\right)  ,O_{i},F_{i}%
:i\in\mathbb{N}\right\}
\]
of complete solutions of the Ricci flow%
\[
\frac{\partial}{\partial t}g_{i}\left(  t\right)  =-2\operatorname*{Rc}\left(
g_{i}\left(  t\right)  \right)
\]
defined for $t\in\left(  \alpha,\omega\right)  $, where $\alpha<0<\omega
\leq\infty$. Each solution is marked by an origin $O_{i}$ and a frame
$F_{i}=\left\{  e_{1}^{i},\dots,e_{n}^{i}\right\}  $ at $O_{i}$ which is
orthonormal with respect to $g_{i}\left(  0\right)  $.

\begin{definition}
We say such a sequence has \textbf{uniformly bounded geometry }if there exist
constants $C_{k}$ for all $k\in\mathbb{N}\cup\left\{  0\right\}  $ such that%
\begin{equation}
\sup_{i\in\mathbb{N}}\sup_{\mathcal{M}_{i}\times\left(  \alpha,\omega\right)
}\left\vert \nabla^{k}\operatorname*{Rm}\left(  g_{i}\right)  \right\vert
_{g_{i}}\leq C_{k}. \label{CurvatureBoundAssumptions}%
\end{equation}

\end{definition}

We denote the eigenvalues of the curvature operator%
\[
\operatorname*{Rm}{}_{i}\left(  x,t\right)  \doteqdot\operatorname*{Rm}\left(
g_{i}\right)  \left(  x,t\right)  :\Lambda^{2}T_{x}\mathcal{M}_{i}%
\rightarrow\Lambda^{2}T_{x}\mathcal{M}_{i}%
\]
by $\lambda_{j}\left(  \operatorname*{Rm}_{i}\right)  \left(  x,t\right)  $,
where $1\leq j\leq m\doteqdot\dim\mathfrak{so}\left(  n\right)  $, and
$\lambda_{1}\leq\cdots\leq\lambda_{m}$.

\begin{definition}
\label{almost nonnegative curvature operators}We call $\left\{  \mathcal{M}%
_{i}^{n},g_{i}\left(  t\right)  ,O_{i},F_{i}:i\in\mathbb{N}\right\}  $ a
\textbf{sequence with almost nonnegative curvature operators} if it has
uniformly bounded geometry and satisfies the following three assumptions:

\begin{enumerate}
\item \label{AlmostNonNegative}there exists a sequence $\delta_{i}\searrow0$
such that
\[
-1\leq-\delta_{i}\leq\lambda_{1}\left(  \operatorname*{Rm}{}_{i}\right)
\left(  x,t\right)
\]
for all $x\in\mathcal{M}_{i}$ and $t\in\left(  \alpha,\omega\right)  $;

\item \label{DiameterGrowth}the manifolds $\left(  \mathcal{M}_{i}^{n}%
,g_{i}\left(  0\right)  \right)  $ are growing without bound,
\[
\lim_{i\rightarrow\infty}\left[  \operatorname*{diam}\left(  \mathcal{M}%
_{i}^{n},g_{i}\left(  0\right)  \right)  \right]  =\infty;
\]
and

\item \label{BumplikePoint}there exists $\varepsilon>0$ such that $O_{i}$ is
an $\varepsilon$-\textbf{bumplike point} at $t=0$, namely
\[
\lambda_{1}\left(  \operatorname*{Rm}{}_{i}\right)  \left(  O_{i},0\right)
\geq\varepsilon.
\]

\end{enumerate}
\end{definition}

\noindent As stated in the introduction, the objective of this paper is to
give a complete proof of:

\begin{theorem}
\label{InjEstForNNCO}For any sequence with almost nonnegative curvature
operators and ~sect$\left(  g_{i}\right)  \left(  x,0\right)  \leq1$ for all
$x\in\mathcal{M}_{i}$ and $i\in\mathbb{N}$, there exists a subsequence
\[
\left\{  \mathcal{M}_{i}^{n},g_{i}\left(  t\right)  ,O_{i},F_{i}\right\}
\]
such that for all $i$,
\[
\operatorname*{inj}{}_{g_{i}\left(  0\right)  }\left(  O_{i}\right)  \geq1.
\]

\end{theorem}

This result is equivalent to Theorem 25.1 of \cite{H1}. The basic strategy of
our proof is the same as the one employed there, and comprises essentially
three steps:

\begin{description}
\item[(a)] use conditions (\ref{AlmostNonNegative}) and (\ref{DiameterGrowth})
to find arbitrarily long minimizing geodesics along which the curvature is
arbitrarily close to nonnegative;

\item[(b)] use condition (\ref{BumplikePoint}) and the strong maximum
principle to construct large uniform neighborhoods of the origins in which the
curvature is uniformly positive; and

\item[(c)] rule out short geodesics in these neighborhoods by means of a
second-variation argument along the long geodesics found in step (a).
\end{description}

However, our implementation of this strategy is distinct in a number of ways
from the methods employed in \cite{H1}. We could not follow one of the steps
in the original proof. In particular, let $\sigma_{i}$ denote the distance to
the cut locus from the origin in $\left(  \mathcal{M}_{i}^{n},g_{i}\left(
0\right)  \right)  $. In the original paper, it is argued that the $\sigma
_{i}$ converge to a continuous function $\ell_{\infty}:\mathcal{S}_{1}%
^{n-1}\rightarrow\left[  0,\infty\right]  $. However there appears to be a gap
in the part of the argument in \cite{H1} that deals with the construction of a
Jacobi field in a geodesic tube for the case that $\exp_{O_{i}}\left(
\ell_{i}V_{i}\right)  =\exp_{O_{i}}\left(  \ell_{i}W_{i}\right)  $ for a
sequence of distinct vectors such that $\left\vert V_{i}-W_{i}\right\vert
\rightarrow0$. This is precisely because collapse for the sequence has not yet
been ruled out at this point in the argument.

\begin{example}
\label{Collapse}Consider a sequence $\left\{  \mathcal{T}_{i}^{2}%
:i=1,2,\dots\right\}  $ of collapsing flat tori with fundamental domains%
\[
\left[  -i,i\right]  \times\left[  -1/i,1/i\right]  \subset\mathbb{R}^{2}.
\]
Take $O_{i}=\left(  0,0\right)  $, and define constant-speed geodesics
\[
\alpha_{i},\beta_{i}:\left[  0,\frac{\sqrt{i^{2}-1}}{i}\right]  \rightarrow
\mathcal{T}_{i}^{2}%
\]
by%
\[
\alpha_{i}\left(  s\right)  =\left(  s,\frac{s}{\sqrt{i^{2}-1}}\right)
\quad\quad\text{and}\quad\quad\beta_{i}\left(  s\right)  =\left(  s,-\frac
{s}{\sqrt{i^{2}-1}}\right)  .
\]
Then $\operatorname{length}\alpha_{i}=\operatorname{length}\beta_{i}=1$ for
all $i$, but their limit in the universal cover $\mathbb{R}^{2}$ is just the
segment $s\mapsto\left(  s,0\right)  $ for $0\leq s\leq1$. Since
$\mathbb{R}^{2}$ is flat, there is no nontrivial Jacobi field which vanishes
at its endpoints.
\end{example}

\begin{remark}
This example does not contain bumplike points, so it is \emph{not }a
counterexample to the claim in \cite{H1}. However, as mentioned above, it does
illustrate difficulties that are due to the possibility of collapse (which is
an issue before an injectivity radius estimate has been proved). One can also
construct `local counterexamples' with constant positive curvature by removing
small neighborhoods of the cone points from $\mathcal{S}^{2}/\mathbb{Z}_{i}$
for $i\in\mathbb{N}$, and letting $i\rightarrow\infty$. This construction does
not produce global counterexamples, since gluing thin infinite cylinders
$\mathcal{S}_{1/\left(  2i\right)  }^{1}\times\left(  0,\infty\right)  $ to
both ends and smoothing the metric will not result in metrics of almost
nonnegative curvature.
\end{remark}

\medskip

In order to overcome this difficulty, we were forced to make some
modifications to both steps (a) and (b) of Hamilton's original proof.

\section{Finding ray-like directions\label{ray-like-directions}}

For each member of the sequence $\left\{  \mathcal{M}_{i}^{n},g_{i}\left(
t\right)  ,O_{i},F_{i}\right\}  $, the frame $F_{i}$ defines a canonical
isometry%
\[
I_{i}:\left(  \mathbb{R}^{n},g_{\operatorname*{can}}\right)  \rightarrow
\left(  T_{O_{i}}\mathcal{M}_{i},g\left(  O_{i},0\right)  \right)  .
\]
Denote the unit sphere bundle of a Riemannian manifold $\left(  \mathcal{M}%
^{n},g\right)  $ by $\mathcal{S}^{n-1}\mathcal{M}^{n}$. For each
$V\in\mathcal{S}_{O_{i}}^{n-1}\mathcal{M}_{i}$, let $\rho_{i}\left(  V\right)
\in(0,\infty]$ denote the distance from $O_{i}$ to the cut point of $O_{i}$
along the geodesic $s\mapsto\exp_{O_{i}}\left(  sV\right)  $ in the metric
$g_{i}\left(  0\right)  $. Denote by $\left(  \mathcal{S}_{1}^{n-1}%
,g_{\operatorname*{can}}\right)  $ the unit sphere in $\mathbb{R}^{n}$ with
its canonical metric, and define%
\[
\sigma_{i}\doteqdot\left.  \rho_{i}\circ I_{i}\right|  _{\mathcal{S}_{1}%
^{n-1}}:\mathcal{S}_{1}^{n-1}\rightarrow(0,\infty].
\]

The set of directions $V\in\mathcal{S}_{1}^{n-1}$ for which $\exp_{O_{i}%
}\left(  sI_{i}\left(  V\right)  \right)  $ is a ray is given by%
\[
\sigma_{i}^{-1}\left(  \infty\right)  =\left\{  V\in\mathcal{S}_{1}^{n-1}%
:\rho_{i}\left(  I_{i}\left(  V\right)  \right)  =\infty\right\}  .
\]
There is no reason to expect $\sigma_{i}^{-1}\left(  \infty\right)  $ to be
nonempty. Indeed, typical applications of Theorem \ref{InjEstForNNCO} are when
$\mathcal{M}_{i}\equiv\mathcal{M}$ is closed or when $\mathcal{M}_{i}$ is
obtained from a closed manifold $\mathcal{M}$ by finitely many surgeries. In
either case we have $\sigma_{i}^{-1}\left(  \infty\right)  \equiv\emptyset$.
Nonetheless, assumption (\ref{DiameterGrowth}) allows us to pick out
directions along which there are arbitrarily long minimizing geodesics.

If $V\in\mathcal{S}_{1}^{n-1}$, let $\mathfrak{S}\left(  V\right)  $ denote
the set of all sequences $\left\{  V_{i}\right\}  \subset\mathcal{S}_{1}%
^{n-1}$ such that $\lim_{i\rightarrow\infty}\left\vert V_{i}-V\right\vert
_{g_{\operatorname*{can}}}=0$. Define
\[
\sigma_{\infty}:\mathcal{S}_{1}^{n-1}\rightarrow\left[  0,\infty\right]
\]
for all $V\in\mathcal{S}_{1}^{n-1}$ by%
\begin{equation}
\fbox{$\sigma_{\infty}\left(  V\right)  \doteqdot\sup_{\mathfrak{S}\left(
V\right)  }\left(  \limsup_{i\rightarrow\infty}\sigma_{i}\left(  V_{i}\right)
\right)  .$} \label{sigma-infinity}%
\end{equation}

\begin{remark}
This definition is the point of departure of our proof from the argument in
\S 25 of \cite{H1}.
\end{remark}

\begin{definition}
If $\left\{  \mathcal{M}_{i}^{n},g_{i}\left(  0\right)  ,O_{i},F_{i}\right\}
$ is a sequence of complete manifolds, its set of \textbf{ray-like directions}
is%
\[
\mathcal{R}_{\infty}\doteqdot\sigma_{\infty}^{-1}\left(  \infty\right)  .
\]

\end{definition}

\noindent In contrast with the sets $\sigma_{i}^{-1}\left(  \infty\right)  $,
the set $\sigma_{\infty}^{-1}\left(  \infty\right)  $ will certainly be nonempty.

\begin{lemma}
\label{R-infinity-is-nonempty}If $\left\{  \mathcal{M}_{i}^{n},g_{i}\left(
0\right)  ,O_{i},F_{i}\right\}  $ is any sequence such that
\[
\operatorname*{diam}\left(  \mathcal{M}_{i}^{n},g_{i}\left(  0\right)
\right)  \rightarrow\infty
\]
as $i\rightarrow\infty$, then $\mathcal{R}_{\infty}$ is nonempty.
\end{lemma}

\begin{proof}
It is a standard fact that each $\sigma_{i}$ is a continuous function on the
compact set $\mathcal{S}_{1}^{n-1}$. For each $i\in\mathbb{N}$, choose
$V_{i}\in\mathcal{S}_{1}^{n-1}$ such that%
\[
\sigma_{i}\left(  V_{i}\right)  =\sup_{V\in\mathcal{S}_{1}^{n-1}}\sigma
_{i}\left(  V\right)  .
\]
Then because $\operatorname*{diam}\left(  \mathcal{M}_{i}^{n},g_{i}\left(
0\right)  \right)  \rightarrow\infty$, we have%
\[
\lim_{i\rightarrow\infty}\sigma_{i}\left(  V_{i}\right)  =\infty.
\]
A subsequence of $V_{i}$ converges to some $V_{\infty}\in\mathcal{S}_{1}%
^{n-1}$. Clearly, $\sigma_{\infty}\left(  V_{\infty}\right)  =\infty$, and
hence $V_{\infty}\in\mathcal{R}_{\infty}$.
\end{proof}

\begin{lemma}
\label{R-infinity-is-compact}If $\left\{  \mathcal{M}_{i}^{n},g_{i}\left(
0\right)  ,O_{i},F_{i}\right\}  $ is any sequence such that
\[
\operatorname*{diam}\left(  \mathcal{M}_{i}^{n},g_{i}\left(  0\right)
\right)  \rightarrow\infty
\]
as $i\rightarrow\infty$, then $\mathcal{R}_{\infty}$ is compact.
\end{lemma}

\begin{proof}
Since $\mathcal{R}_{\infty}\subseteq\mathcal{S}_{1}^{n-1}$, it will suffice to
show that $\mathcal{R}_{\infty}$ is closed. So let $\left\{  V^{\alpha}%
:\alpha\in\mathbb{N}\right\}  $ be a sequence from $\mathcal{R}_{\infty}$ such
that $\lim_{\alpha\rightarrow\infty}V^{\alpha}=V^{\infty}\in\mathcal{S}%
_{1}^{n-1}$ exists. Then by definition of $\mathcal{R}_{\infty}$, there is for
every $\alpha$ a sequence $\left\{  V_{i}^{\alpha}:i\in\mathbb{N}\right\}  $
of unit vectors such that%
\[
\lim_{i\rightarrow\infty}V_{i}^{\alpha}=V^{\alpha}\quad\quad\text{and}%
\quad\quad\lim_{i\rightarrow\infty}\sigma_{i}\left(  V_{i}^{\alpha}\right)
=\infty.
\]
Observe that we can choose $i\left(  \alpha\right)  \geq\alpha$ for all
$\alpha$ such that%
\[
\left|  V_{i\left(  \alpha\right)  }^{\alpha}-V_{\alpha}\right|
_{g_{\operatorname*{can}}}<\frac{1}{\alpha}\quad\quad\text{and}\quad
\quad\sigma_{i\left(  \alpha\right)  }\left(  V_{i\left(  \alpha\right)
}^{\alpha}\right)  >\alpha.
\]
Then%
\[
\lim_{\alpha\rightarrow\infty}V_{i\left(  \alpha\right)  }^{a}=\lim
_{\alpha\rightarrow\infty}V^{\alpha}=V^{\infty}%
\]
and%
\[
\lim_{\alpha\rightarrow\infty}\sigma_{i\left(  \alpha\right)  }\left(
V_{i\left(  \alpha\right)  }^{\alpha}\right)  =\infty.
\]
Hence $V^{\infty}\in\mathcal{R}_{\infty}$.
\end{proof}

\begin{remark}
The reader may wish to contrast $\mathcal{R}_{\infty}$ with the set
$\mathcal{D}$ of \emph{distinguished directions }defined in \cite{H1} as%
\[
\mathcal{D}\doteqdot\ell_{\infty}^{-1}\left(  \infty\right)  ,
\]
where $\ell_{\infty}:\mathcal{S}_{1}^{n-1}\rightarrow\left[  0,\infty\right]
$ is defined by%
\[
\ell_{\infty}\left(  V\right)  \doteqdot\lim_{i\rightarrow\infty}\sigma
_{i}\left(  V_{i}\right)  .
\]
It is claimed there that this limit is independent of the sequence
$V_{i}\rightarrow V$. In contrast, our proof does not rely on such a property
of independence of sequence.
\end{remark}

\begin{remark}
[Tilting frames]\label{TiltingFrames}Without the bump-like origins condition
in Definition \ref{almost nonnegative curvature operators}, it is not
necessarily true that the distance with respect to $g_{i}\left(  0\right)  $
from $O_{i}$ to the cut locus in a direction $V\in\mathcal{R}_{\infty}$ can be
made arbitrarily large by going far enough out in the sequence. This is
illustrated by the following example: for $i=1,2,\dots$, let
\[
\left(  \mathcal{M}_{i}^{2},g_{i}\right)  =\mathcal{S}_{1/i}^{1}%
\times\mathbb{R},
\]
where
\[
\mathcal{S}_{1/i}^{1}=\left\{  \left(  x,y\right)  \in\mathbb{R}^{2}%
:x^{2}+y^{2}=1/i^{2}\right\}  .
\]
It is not important that the manifolds $\mathcal{M}_{i}$ are not compact,
since once can replace the $\mathcal{M}_{i}$ by tori of various lengths.
Regard each $\mathcal{M}_{i}$ as a submanifold of $\mathbb{R}^{3}$ with
coordinates $\left(  x,y,z\right)  $, and take the origins to be
$O_{i}=\left(  0,1/i,0\right)  $. (Actually, any sequence of points in
$\mathcal{M}_{i}$ will do.) Given any positive real number $\lambda$, define
the tilted frame $\mathcal{F}_{i}=\left\{  e_{1}^{i},e_{2}^{i}\right\}  $ at
$O_{i}$ by rotating the standard frame
\[
\mathcal{F}=\left\{  e_{1}=\left(  0,0,1\right)  ,\,e_{2}=\left(
1,0,0\right)  \right\}
\]
clockwise by an angle of $\arctan\left(  \pi/i\lambda\right)  $. Using $I_{i}$
to identify $\mathbb{R}^{2}$ with $T_{O_{i}}\mathcal{M}_{i}$, we have
$I_{i}^{-1}\left(  e_{1}^{i}\right)  \equiv\left(  1,0\right)  \in
\mathbb{R}^{2}$. Hence
\[
\sigma_{i}\left(  \left(  1,0\right)  \right)  =\sigma_{i}\left(  I_{i}%
^{-1}\left(  e_{1}^{i}\right)  \right)  =\sqrt{\lambda^{2}+\pi^{2}/i^{2}},
\]
so that
\[
\lim_{i\rightarrow\infty}\sigma_{i}\left(  \left(  1,0\right)  \right)
=\lambda.
\]
On the other hand%
\[
\sigma_{i}\left(  I_{i}^{-1}\left(  e_{1}\right)  \right)  =\infty.
\]
Since $\lim_{i\rightarrow\infty}I_{i}^{-1}\left(  e_{1}\right)  =\left(
1,0\right)  $, this implies that%
\[
\sigma_{\infty}\left(  \left(  1,0\right)  \right)  =\infty.
\]
Observe in particular that%
\[
\lim_{i\rightarrow\infty}\sigma_{i}\left(  \left(  1,0\right)  \right)
\neq\sigma_{\infty}\left(  \left(  1,0\right)  \right)  .
\]

\end{remark}%

\begin{center}
\includegraphics[
trim=0.000000in 0.000000in 0.742140in 8.461139in,
height=1.3647in,
width=4.1554in
]%
{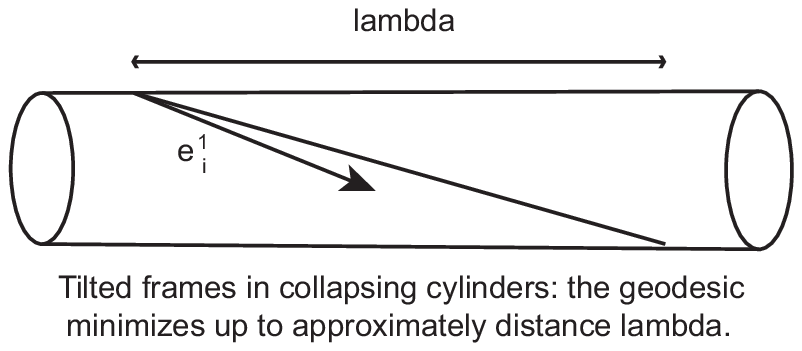}%
\end{center}
%

\begin{center}
\includegraphics[
trim=0.000000in 0.000000in 0.742140in 8.461139in,
height=1.3647in,
width=4.1554in
]%
{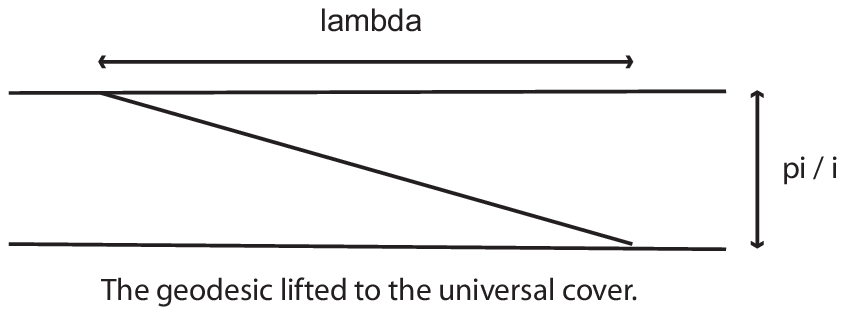}%
\end{center}

\begin{remark}
If $V\in\mathcal{R}_{\infty}$, it \emph{is} true that we can find arbitrarily
long minimizing geodesics in directions arbitrarily close to $V$. In
particular, $V\in\mathcal{S}_{1}^{n-1}$ belongs to $\mathcal{R}_{\infty}$ if
and only if there is a sequence $\left\{  V_{i}\right\}  $ from $\mathcal{S}%
_{1}^{n-1}$ such that $\lim_{i\rightarrow\infty}\left\vert V_{i}-V\right\vert
_{g_{\operatorname*{can}}}=0$ and $\lim_{i\rightarrow\infty}\sigma_{i}\left(
V_{i}\right)  =\infty$.
\end{remark}

\begin{remark}
\label{R-infinity-unstable}\emph{A priori, }our definition allows
$\mathcal{R}_{\infty}$ to become smaller each time we pass to a subsequence.
We shall deal with this issue carefully in the sequel.
\end{remark}

\section{Finding large balls of positive curvature\label{large balls}}

Our proof of Theorem \ref{InjEstForNNCO} requires us to show that the
curvature can be made positive in arbitrarily large neighborhoods of the
origin by going sufficiently far out in the sequence. The key to this result
is the strong maximum principle of \cite{H6}, which lets us pass from purely
local results to results that hold on arbitrarily large sets. Because this
part of our argument is essentially the same as the construction in \S 25 of
\cite{H1}, we shall omit or merely outline most proofs.

\subsection{Preconvergence in geodesic tubes}

To begin, we recall the procedure of taking limits of the pullback metrics in
geodesic tubes, where one can avoid the need for an injectivity radius
estimate for the sequence $\left\{  \mathcal{M}_{i}^{n},g_{i}\left(  t\right)
,O_{i},F_{i}\right\}  $.

If $\gamma:\left(  -L,L\right)  \rightarrow\left(  \mathcal{M}^{n},g\right)  $
is a unit-speed geodesic, we denote its normal bundle by%
\[
N_{\gamma}\doteqdot\left\{  V\in\gamma^{\ast}\left(  T\mathcal{M}\right)
:\left\langle V,\dot{\gamma}\right\rangle =0\right\}  .
\]
Note that $N_{\gamma}$ is a rank $n-1$ vector bundle over $\left(
-L,L\right)  $. Define $\Phi:N_{\gamma}\rightarrow\mathcal{M}$ by%
\[
\Phi\left(  \gamma\left(  s\right)  \right)  \doteqdot\left.  \exp
_{\gamma\left(  s\right)  }\right|  _{N_{\gamma}\left(  \gamma\left(
s\right)  \right)  }.
\]
Now let $F=\left\{  e_{1},\dots,e_{n}\right\}  $ be any orthonormal frame at
$\gamma\left(  0\right)  $ with $e_{1}=\dot{\gamma}$. Taking the pullback of
$F$ and parallel translating it along $\gamma$, we obtain a orthonormal basis
in each fiber of $N_{\gamma}$, which we continue to denote by $\left\{
e_{2},\dots,e_{n}\right\}  $. Denoting the open ball of radius $r>0$ centered
at $\vec{0}\in\mathbb{R}^{n-1}$ by $B(\vec{0},r)$, we define the cylinder
$T_{L,r}\doteqdot\left(  -L,L\right)  \times B(\vec{0},r)$ and a map
$\iota_{F}:T_{L,r}\rightarrow N_{\gamma}$ by%
\[
\iota_{F}:\left(  v^{1},\;v^{2},\dots,v^{n}\right)  \mapsto\left(
\gamma\left(  v^{1}\right)  ,%
{\textstyle\sum_{k=2}^{n}}
v^{k}e_{k}\left(  v_{1}\right)  \right)  ,
\]
where $\left(  v^{1},\;v^{2},\dots,v^{n}\right)  $ are the natural Euclidean
coordinates in $T_{L,r}$.

The following two results are well known; their proofs are virtually identical
to those of corresponding results for exponential coordinates.

\begin{lemma}
If $-1\leq\operatorname{sect}\left(  g\right)  \leq1$, there exists $\rho>0$
depending only on $n$ such that the map
\[
\Psi_{F,L,2\rho}\doteqdot\Phi\circ\iota_{F}:T_{L,2\rho}\rightarrow
\mathcal{M}^{n}%
\]
is an immersion for every geodesic $\gamma:\left(  -L,L\right)  \rightarrow
\mathcal{M}$, every orthonormal frame $F=\left\{  \dot{\gamma}=e_{1}%
,e_{2},\dots,e_{n}\right\}  $ at $\gamma\left(  0\right)  $, and every $L>0$.
\end{lemma}

\begin{lemma}
\label{EstForGeodesicTube}Suppose $-1\leq\operatorname{sect}\left(  g\right)
\leq1$. If $\delta$ denotes the Euclidean metric on $T_{L,2\rho}$ and $h$ is
the pullback metric%
\[
h\doteqdot\Psi_{F,L,2\rho}^{\ast}g,
\]
then there exist constants $0<c<C<\infty$ depending only on $n$ such that%
\begin{equation}
c\,\delta\leq h\leq C\,\delta. \label{MetricEquivalence}%
\end{equation}
Moreover, if there are constants $C_{k}$ such that%
\[
\left|  \nabla^{k}\operatorname*{Rm}\right|  _{g}\leq C_{k}%
\]
for all $k\in\mathbb{N}$, then there exist $C_{k}^{\prime}=C_{k}^{\prime
}\left(  C_{0},\dots,C_{k}\right)  $ such that%
\begin{equation}
\left|  \frac{\partial^{k}}{\partial v^{\alpha_{1}}\cdots\partial
v^{\alpha_{k}}}h\right|  _{\delta}\leq C_{k}^{\prime}.
\label{DerivativeEstimates}%
\end{equation}

\end{lemma}

Now let $\left\{  \mathcal{M}_{i}^{n},g_{i}\left(  t\right)  ,O_{i}%
,F_{i}\right\}  $ be a sequence having uniformly bounded geometry. For each
$A\in\operatorname*{O}\left(  n\right)  $ and $L>0$, there is a sequence
$\left(  T_{L,2\rho},\;\Psi_{AF_{i},L,2\rho}^{\ast}\left(  g_{i}\left(
t\right)  \right)  \right)  $ of geodesic tubes extending in the direction
determined by the frame $AF_{i}$ obtained from the natural action of $A$ on
$F_{i}$. Let $\left\{  A_{\alpha}:\alpha\in\mathbb{N}\right\}  $ be a
countable dense set in the compact Lie group $\operatorname*{O}\left(
n\right)  $, and let $\left\{  L_{\beta}:\beta\in\mathbb{N}\right\}  $ be a
sequence of positive real numbers with $L_{\beta}\nearrow\infty$ as
$\beta\rightarrow\infty$.

\begin{lemma}
There exists a subsequence $\left\{  \mathcal{M}_{i}^{n},g_{i}\left(
t\right)  ,O_{i},F_{i}\right\}  $ such that for all $\alpha,\beta\in
\mathbb{N}$, the geodesic tube
\[
\left(  T_{L_{\beta},2\rho},\;\Psi_{A_{\alpha}F_{i},L_{\beta},2\rho}^{\ast
}\left(  g_{i}\left(  t\right)  \right)  \right)
\]
converges as $i\rightarrow\infty$ uniformly in each $C^{k}$ norm to a
solution
\[
\left(  T_{L_{\beta},2\rho},\;h_{A_{\alpha},L_{\beta},2\rho}^{\infty}\left(
t\right)  \right)
\]
of the Ricci flow for $t\in\left(  \alpha,\omega\right)  $.
\end{lemma}

\begin{proof}
By Lemma \ref{EstForGeodesicTube}, $\Psi_{A_{\alpha}F_{i},L_{\beta},\rho
}^{\ast}\left(  g_{i}\left(  t\right)  \right)  $ is a sequence of solutions
of the Ricci flow on $T_{L_{\beta},\rho}$ such that the pullback metrics
satisfy (\ref{MetricEquivalence}) and (\ref{DerivativeEstimates}) uniformly in
$i\in\mathbb{N}$ and $t\in\left(  \alpha,\omega\right)  $. Thus the result
follows from Arzela--Ascoli by consecutive diagonalization arguments.
\end{proof}

Since $L_{\beta}\rightarrow\infty$, we have as an immediate consequence:

\begin{corollary}
There exists a subsequence $\left\{  \mathcal{M}_{i}^{n},g_{i}\left(
t\right)  ,O_{i},F_{i}\right\}  $ such that for all $\alpha\in\mathbb{N}$ and
all $L>0$, the geodesic tube
\[
\left(  T_{L,2\rho},\;\Psi_{A_{\alpha}F_{i},L,2\rho}^{\ast}\left(
g_{i}\left(  t\right)  \right)  \right)
\]
converges as $i\rightarrow\infty$ uniformly in each $C^{k}$ norm to a
solution
\[
\left(  T_{L,2\rho},\;h_{A_{\alpha},L,2\rho}^{\infty}\left(  t\right)
\right)
\]
of the Ricci flow for $t\in\left(  \alpha,\omega\right)  $.
\end{corollary}

From this, one can obtain a subsequence that converges in any geodesic tube of
any finite length.

\begin{definition}
A sequence $\left\{  \mathcal{M}_{i}^{n},g_{i}\left(  t\right)  ,O_{i}%
,F_{i}\right\}  $ that has\ uniformly bounded geometry is said to be
\textbf{preconvergent in geodesic tubes }if for all $A\in\operatorname*{O}%
\left(  n\right)  $ and $L>0$, the geodesic tube%
\[
\left(  T_{L,\rho},\;\Psi_{AF_{i},L,\rho}^{\ast}\left(  g_{i}\left(  t\right)
\right)  \right)
\]
converges as $i\rightarrow\infty$ uniformly in each $C^{k}$ norm to a solution%
\[
\left(  T_{L,\rho},\;h_{A,L,\rho}^{\infty}\left(  t\right)  \right)
\]
of the Ricci flow for $t\in\left(  \alpha,\omega\right)  $.
\end{definition}

\begin{lemma}
\label{ConvergenceInGeodesicTubes}Any sequence having uniformly bounded
geometry contains a subsequence that is preconvergent in geodesic tubes.
\end{lemma}

\begin{proof}
Given $A\in\operatorname*{O}\left(  n\right)  $ and $L>0$, consider
$\Psi_{AF_{i},L,\rho}:T_{L,\rho}\rightarrow\left(  \mathcal{M}_{i}%
,g_{i}\left(  t\right)  \right)  $. Choose a sequence $A_{\alpha}$ such that
$A_{\alpha}\rightarrow A$ as $\alpha\rightarrow\infty$ in some (hence any)
norm on $\operatorname*{O}\left(  n\right)  $. By standard covering-space
theory, there exists $\alpha^{\prime}=\alpha^{\prime}\left(  L,\rho\right)  $
independent of $i$ such that for all $\alpha\geq\alpha^{\prime}$ there exists
a smooth embedding $\iota_{\alpha,L,\rho}:T_{L,\rho}\rightarrow T_{2L,2\rho}$
such that $\Psi_{A_{\alpha}F_{i},2L,2\rho}\circ\iota_{\alpha,L,\rho}%
=\Psi_{AF_{i},L,\rho}$. Note that all derivatives of $\iota_{\alpha,L,\rho}$
are bounded uniformly with respect to $\alpha$. Note too that as
$\alpha\rightarrow\infty$, the maps $\iota_{\alpha,L,\rho}$ converge uniformly
in each $C^{k}$ norm to the inclusion map $\iota_{L,\rho}:T_{L,\rho
}\rightarrow T_{2L,2\rho}$ that\ is defined so that $\Psi_{AF_{i},2L,2\rho
}\circ\iota_{L,\rho}=\Psi_{AF_{i},L,\rho}$. So as $\alpha\rightarrow\infty$,
we have $\left(  T_{L,\rho},\;\Psi_{A_{\alpha}F_{i},L,\rho}^{\ast}\left(
g_{i}\left(  t\right)  \right)  \right)  \rightarrow\left(  T_{L,\rho}%
,\;\Psi_{AF_{i},L,\rho}^{\ast}\left(  g_{i}\left(  t\right)  \right)  \right)
$ uniformly in each $C^{k}$ norm. Since this convergence is independent of
$i$, a routine diagonalization argument completes the proof.
\end{proof}

\begin{remark}
The solutions $\left(  T_{L,\rho},\;h_{A,L,\rho}^{\infty}\left(  t\right)
\right)  $ obtained by this construction are not complete.
\end{remark}

\subsection{Preconvergence in distance}

The distance function in each geodesic tube gives an upper bound for the
distance function in the original geometry. This fact ensures that any
sequence that is preconvergent in geodesic tubes contains a subsequence that
is \emph{preconvergent in distance, }namely a subsequence such that the limit
\[
d_{\infty}\left(  X,Y\right)  \doteqdot\lim_{i\rightarrow\infty}\left[
d_{i}\left(  \exp_{O_{i}}\left(  I_{i}\left(  X\right)  \right)  ,\exp_{O_{i}%
}\left(  I_{i}\left(  Y\right)  \right)  \right)  \right]
\]
exists for all $X,Y\in\mathbb{R}$, where $d_{i}$ denotes the distance function
induced on the manifold $\mathcal{M}_{i}^{n}$ by the Riemannian metric
$g_{i}(0)$. Preconvergence in distance is a stability property; it ensures,
for instance, that the distance to the cut locus in a particular direction is
not going to infinity along one subsequence while remaining uniformly bounded
along another subsequence. However, the methods we develop in
\S \ \ref{boundedness} --- in particular, covering $\mathcal{R}_{\infty}$ by
a\ finite $\varepsilon$-net of directions for which the $\limsup$ in
Definition (\ref{sigma-infinity}) is attained as a limit --- make it
unnecessary to use this property. Consequently, we omit the proof.

\subsection{Preconvergence to positive curvature}

\begin{notation}
If $x\in\mathcal{M}_{i}^{n}$, we denote the open $g_{i}\left(  0\right)
$-ball of radius $r>0$ centered at $x$ by%
\[
B_{i}\left(  x,r\right)  \doteqdot\left\{  y\in\mathcal{M}_{i}^{n}%
:d_{i}\left(  x,y\right)  <r\right\}
\]
and the corresponding closed ball by%
\[
\bar{B}_{i}\left(  x,r\right)  \doteqdot\left\{  y\in\mathcal{M}_{i}^{n}%
:d_{i}\left(  x,y\right)  \leq r\right\}  .
\]
We denote by $d_{i}$ the distance function induced on $\mathcal{M}_{i}^{n}$ by
the metric $g_{i}\left(  0\right)  $.
\end{notation}

\begin{definition}
We say a sequence $\left\{  \mathcal{M}_{i}^{n},g_{i}\left(  t\right)
,O_{i},F_{i}\right\}  $ with almost nonnegative curvature operators is
\textbf{preconverging to positive curvature} if for each $L>0$ there are some
$\eta\left(  L\right)  >0$ and $\iota\left(  L\right)  \in\mathbb{N}$ such
that%
\[
\lambda_{1}\left(  \operatorname*{Rm}{}_{i}\right)  \left(  x,0\right)
\geq\eta
\]
for all $i\geq\iota$ and all $x\in\bar{B}_{i}\left(  O_{i},L\right)  $.
\end{definition}

\begin{proposition}
\label{CurvatureBoundFromBelow}Any sequence with almost nonnegative curvature
operators contains a subsequence that is preconverging to positive curvature.
\end{proposition}

\begin{proof}
Since the result is equivalent to Lemma 25.2 of \cite{H1}, we shall merely
sketch the proof. By Lemma \ref{ConvergenceInGeodesicTubes}, we first pass to
a subsequence that is preconvergent in geodesic tubes. If the proposition is
false, there is some $L\in\left(  0,\infty\right)  $ such that for all
$\eta>0$ and each $i_{0}\in\mathbb{N}$, we have $\lambda_{1}\left(
\operatorname*{Rm}{}_{i}\right)  \left(  x,0\right)  <\eta$ for some $i\geq
i_{0}$ and some $x\in\bar{B}_{i}\left(  O_{i},L\right)  $. One can then argue
that there exists a further subsequence such that
\[
\lambda_{1}\left(  \operatorname*{Rm}{}_{i}\right)  \left(  \exp_{O_{i}%
}\left(  V_{i}\right)  ,0\right)  \rightarrow0
\]
as $i\rightarrow\infty$, where $I_{i}^{-1}\left(  V_{i}\right)  $ converges to
some $V_{\infty}\in\mathbb{R}^{n}$ with $0<\left|  V_{\infty}\right|  \leq L$.
Now preconvergence in geodesic tubes ensures that the subsequence to which we
have passed has the property that%
\[
\left(  T_{2L,\rho},\;\Psi_{AF_{i},2L,\rho}^{\ast}\left(  g_{i}\left(
t\right)  \right)  \right)  \rightarrow\left(  T_{2L,\rho},\;h_{A,2L,\rho
}^{\infty}\left(  t\right)  \right)  ,
\]
where $A\in\operatorname*{O}\left(  n\right)  $ is chosen such that the
geodesic tube lies in the direction $V_{\infty}/\left|  V_{\infty}\right|  $.
It follows that
\[
\lambda_{1}\left(  \operatorname*{Rm}\left(  h_{A,2L,\rho}^{\infty}\left(
0\right)  \right)  \right)  \left(  \left|  V_{\infty}\right|  ,\vec
{0}\right)  =0.
\]
Because assumption (\ref{AlmostNonNegative}) implies that $\operatorname*{Rm}%
\left(  h_{A,2L,\rho}^{\infty}\left(  0\right)  \right)  \geq0$, we may then
apply the strong maximum principle in the form proved in \cite{H6} to conclude
that
\[
\lambda_{1}\left(  \operatorname*{Rm}\left(  h_{A,2L,\rho}^{\infty}\left(
0\right)  \right)  \right)  \left(  0,\vec{0}\right)  =0.
\]
But this is possible only if $\lambda_{1}\left(  \operatorname*{Rm}\left(
\Psi_{AF_{i},2L,\rho}^{\ast}\left(  g_{i}\left(  0\right)  \right)  \right)
\right)  \left(  0,\vec{0}\right)  \rightarrow0$ as $i\rightarrow\infty$,
which contradicts assumption (\ref{BumplikePoint}), because%
\[
\lambda_{1}\left(  \operatorname*{Rm}\left(  \Psi_{AF_{i},2L,\rho}^{\ast
}\left(  g_{i}\left(  0\right)  \right)  \right)  \right)  \left(  0,\vec
{0}\right)  =\lambda_{1}\left(  \operatorname*{Rm}{}_{i}\right)  \left(
O_{i},0\right)  \geq\varepsilon>0.
\]
The contradiction proves the proposition.
\end{proof}

\section{Mimicking the sublevel sets of a Busemann function\label{mimicking}}

If $\left(  \mathcal{M}^{n},g\right)  $ is a complete noncompact manifold of
positive sectional curvature bounded above by $\kappa<\infty$, then Gromoll
and Meyer \cite{GM} proved that its injectivity radius can be bounded below by
$\pi/\sqrt{\kappa}$. One way to prove this is to fix an origin $O\in
\mathcal{M}^{n}$, use the rays emanating from $O$ to construct a Busemann
function associated to that origin, use that Busemann function to construct a
totally convex neighborhood $N$ of $O$, and then use a second variation
argument along rays to rule out short geodesics in the neighborhood $N$. (See
Greene \cite{G} for a survey of noncompact manifolds with nonnegative curvature.)

Following Hamilton, we want to mimic this construction along a sequence
$\left\{  \mathcal{M}_{i}^{n},g_{i}\left(  t\right)  ,O_{i},F_{i}\right\}  $
that is preconverging to positive curvature. Since we may only have ray-like
directions in general, we need a substitute for the Busemann construction.

\begin{notation}
We shall henceforth identify $V\in\mathbb{R}^{n}$ with $I_{i}\left(  V\right)
\in T_{O_{i}}\mathcal{M}_{i}$ and vice versa, using the canonical isometries
\[
I_{i}:\left(  \mathbb{R}^{n},g_{\operatorname*{can}}\right)  \rightarrow
\left(  T_{O_{i}}\mathcal{M}_{i},g\left(  O_{i},0\right)  \right)  .
\]

\end{notation}

As a substitute for the sublevel sets of a Busemann function, we define%

\begin{align}
&  N_{i}\left(  L,K\right) \label{Busemann}\\
&  \doteqdot\left\{  \exp_{O_{i}}\left(  tW\right)  \left|
\begin{array}
[c]{c}%
W\in\mathcal{S}_{1}^{n-1}\text{ and }t\in\left[  0,\sigma_{i}\left(  W\right)
\right]  \backslash\left\{  \infty\right\}  \text{ are such that}\\
\mathstrut\\
\text{for every }V\in\mathcal{S}_{1}^{n-1}\text{ with }\sigma_{i}\left(
V\right)  \geq L\text{,}\\
\mathstrut\\
\text{all }r\in\left[  0,\sigma_{i}\left(  V\right)  \right]  \backslash
\left\{  \infty\right\}  \text{, and all }s\in\left[  0,t\right]  \text{, we
have}\\
\mathstrut\\
d_{i}\left(  \exp_{O_{i}}\left(  rV\right)  ,\exp_{O_{i}}\left(  sW\right)
\right)  \geq r-K
\end{array}
\right.  \right\}  .\nonumber
\end{align}
The corresponding sets in \S 25 of \cite{H1} are constructed using those $V$
for which $\ell_{\infty}\left(  V\right)  =\infty$. But for our proof, it is
important to allow $V$ such that $\sigma_{i}\left(  V\right)  $ is large but finite.

Notice\ that each $N_{i}\left(  L,K\right)  $ is weakly star shaped with
respect to $O_{i}$: namely, $\exp_{O_{i}}\left(  tW\right)  \in N_{i}\left(
L,K\right)  $ implies that $\exp_{O_{i}}\left(  sW\right)  \in N_{i}\left(
L,K\right)  $ for all $0\leq s\leq t$. It will be useful to collect a few more
elementary observations about the sets $N_{i}\left(  L,K\right)  $.

\begin{lemma}
\label{N_iLK-contains-K-ball}$N_{i}\left(  L,K\right)  $ contains the closed
ball of radius $\min\left\{  \pi,K\right\}  $:%
\[
\bar{B}_{i}\left(  O_{i},\min\left\{  \pi,K\right\}  \right)  \subseteq
N_{i}\left(  L,K\right)  .
\]

\end{lemma}

\begin{proof}
By assumption (\ref{AlmostNonNegative}), there are no conjugate points of
$O_{i}$ in $B_{i}\left(  O_{i},\pi\right)  $. So if $x\in\bar{B}_{i}\left(
O_{i},\min\left\{  \pi,K\right\}  \right)  $, there exist $W\in\mathcal{S}%
_{1}^{n-1}$ and $t\in\left[  0,\sigma_{i}\left(  W\right)  \right]  $ such
that $x=\exp_{O_{i}}\left(  tW\right)  $. Then for any $V\in\mathcal{S}%
_{1}^{n-1}$ with $\ \sigma_{i}\left(  V\right)  \geq L$, all $r\in\left[
0,\sigma_{i}\left(  V\right)  \right]  \backslash\left\{  \infty\right\}  $,
and all $s\in\left[  0,t\right]  $, the triangle inequality yields%
\begin{align*}
r  &  =d_{i}\left(  O_{i},\exp_{O_{i}}\left(  rV\right)  \right) \\
&  \leq d_{i}\left(  O_{i},\exp_{O_{i}}\left(  sW\right)  \right)
+d_{i}\left(  \exp_{O_{i}}\left(  sW\right)  ,\exp_{O_{i}}\left(  rV\right)
\right) \\
&  =s+d_{i}\left(  \exp_{O_{i}}\left(  sW\right)  ,\exp_{O_{i}}\left(
rV\right)  \right) \\
&  \leq K+d_{i}\left(  \exp_{O_{i}}\left(  sW\right)  ,\exp_{O_{i}}\left(
rV\right)  \right)  .
\end{align*}
The inequality on the last line holds because $s\leq t=d_{i}\left(
O_{i},x\right)  \leq\min\left\{  \pi,K\right\}  $.
\end{proof}

\begin{lemma}
\label{Tangential-Radial}If $\exp_{O_{i}}\left(  tW\right)  \in N_{i}\left(
L,K\right)  $, then for all $V\in\mathcal{S}_{1}^{n-1}$ such that $\sigma
_{i}\left(  V\right)  \geq L$, all $r\in\left[  0,\sigma_{i}\left(  V\right)
\right]  \backslash\left\{  \infty\right\}  $, and all $s\in\left[
0,t\right]  $, we have%
\[
d_{i}\left(  O_{i},\exp_{O_{i}}\left(  sW\right)  \right)  \leq K+2\cdot
d_{i}\left(  \exp_{O_{i}}\left(  sW\right)  ,\exp_{O_{i}}\left(  rV\right)
\right)  .
\]

\end{lemma}

\begin{proof}
By the triangle inequality,%
\[
d_{i}\left(  O_{i},\exp_{O_{i}}\left(  rV\right)  \right)  \geq d_{i}\left(
O_{i},\exp_{O_{i}}\left(  sW\right)  \right)  -d_{i}\left(  \exp_{O_{i}%
}\left(  rV\right)  ,\exp_{O_{i}}\left(  sW\right)  \right)  .
\]
So because $\exp_{O_{i}}\left(  tW\right)  \in N_{i}\left(  L,K\right)  $, we
have%
\begin{align*}
&  d_{i}\left(  O_{i},\exp_{O_{i}}sW\right)  -d_{i}\left(  \exp_{O_{i}}\left(
rV\right)  ,\exp_{O_{i}}\left(  sW\right)  \right) \\
&  \leq d_{i}\left(  O_{i},\exp_{O_{i}}\left(  rV\right)  \right)  =r\\
&  \leq K+d_{i}\left(  \exp_{O_{i}}\left(  rV\right)  ,\exp_{O_{i}}\left(
sW\right)  \right)  ,
\end{align*}
and hence $d_{i}\left(  O_{i},\exp_{O_{i}}\left(  sW\right)  \right)  \leq
K+2\cdot d_{i}\left(  \exp_{O_{i}}\left(  rV\right)  ,\exp_{O_{i}}\left(
sW\right)  \right)  $.
\end{proof}

\begin{corollary}
\label{LK-ball-estimate}$N_{i}\left(  L,K\right)  $ is contained in the closed
ball of radius $\max\left\{  L,K\right\}  $:
\[
N_{i}\left(  L,K\right)  \subseteq\bar{B}_{i}\left(  O_{i},\max\left\{
L,K\right\}  \right)  .
\]

\end{corollary}

\begin{proof}
If $\exp_{O_{i}}\left(  tW\right)  \in N_{i}\left(  L,K\right)  $ and
$d_{i}\left(  O_{i},\exp_{O_{i}}\left(  tW\right)  \right)  >L$, apply Lemma
\ref{Tangential-Radial} with $V=W$ and $r=s=t$ to get $d_{i}\left(  O_{i}%
,\exp_{O_{i}}\left(  tW\right)  \right)  \leq K$.
\end{proof}

\begin{lemma}
\label{N_iLK-are-compact}Each $N_{i}\left(  L,K\right)  $ is compact.
\end{lemma}

\begin{proof}
By Corollary \ref{LK-ball-estimate}, it suffices to show that each
$N_{i}\left(  L,K\right)  $ is closed. Let
\[
\left\{  x_{\alpha}=\exp_{O_{i}}\left(  t_{\alpha}W_{\alpha}\right)
:\alpha\in\mathbb{N}\right\}  \subset N_{i}\left(  L,K\right)
\]
be a sequence such that $\lim_{\alpha\rightarrow\infty}x_{\alpha}\doteqdot
x\in\mathcal{M}_{i}^{n}$ exists. Then%
\[
\lim_{\alpha\rightarrow\infty}t_{\alpha}\doteqdot t\in\left[  0,\max\left\{
L,K\right\}  \right]
\]
and%
\[
\lim_{\alpha\rightarrow\infty}W_{\alpha}\doteqdot W\in\mathcal{S}_{1}^{n-1}%
\]
exist also. Since $t_{\alpha}\in\left[  0,\sigma_{i}\left(  W_{\alpha}\right)
\right]  $ for all $\alpha$, the continuity of $\sigma_{i}$ implies that
$t\in\left[  0,\sigma_{i}\left(  W\right)  \right]  $.

Now let $V\in\mathcal{S}_{1}^{n-1}$ with $\sigma_{i}\left(  V\right)  >L$,
$r\in\left[  0,\sigma_{i}\left(  V\right)  \right]  \backslash\left\{
\infty\right\}  $, and $s\in\left[  0,t\right]  $ be given. Choose $s_{\alpha
}\in\left[  0,t_{\alpha}\right]  $ such that $\lim_{\alpha\rightarrow\infty
}s_{\alpha}=s$. Then for any $\varepsilon>0$, there exists $A<\infty$ such
that for all $\alpha\geq A$ we have%
\[
\left|  d_{i}\left(  \exp_{O_{i}}\left(  rV\right)  ,\exp_{O_{i}}\left(
sW\right)  \right)  -d_{i}\left(  \exp_{O_{i}}\left(  rV\right)  ,\exp_{O_{i}%
}\left(  s_{\alpha}W_{\alpha}\right)  \right)  \right|  <\varepsilon.
\]
Since $\exp_{O_{i}}\left(  t_{\alpha}W_{\alpha}\right)  \in N_{i}\left(
L,K\right)  $ for each $\alpha$, this implies that%
\[
d_{i}\left(  \exp_{O_{i}}\left(  rV\right)  ,\exp_{O_{i}}\left(  sW\right)
\right)  \geq r-K-\varepsilon.
\]
Since $\varepsilon>0$ was arbitrary, we see that $x=\exp_{O_{i}}\left(
tW\right)  \in N_{i}\left(  L,K\right)  $.
\end{proof}

\section{The boundedness property\label{boundedness}}

What is not obvious is the fact that the $N_{i}\left(  L,K\right)  $ can be
uniformly bounded. That is the content of the following crucial result:

\begin{proposition}
[Boundedness property]\label{N_iLK-are-bounded}Any sequence preconverging to
positive curvature contains a subsequence for which there exists a constant
$C<\infty$ depending on $K$ such that for each $L\in\left(  0,\infty\right)
$, there exists $I\left(  L\right)  $ such that for all $i\geq I\left(
L\right)  $, we have%
\[
N_{i}\left(  L,K\right)  \subseteq B_{i}\left(  O_{i},C\right)  .
\]

\end{proposition}

The proof of the proposition has two main steps. The first step is to observe
that by passing to a subsequence, we can in a sense replace $\mathcal{R}%
_{\infty}$ by a finite $\varepsilon$-net of directions for which the $\limsup$
in Definition (\ref{sigma-infinity}) is attained as a limit. We have already
observed in Remark \ref{R-infinity-unstable} that $\mathcal{R}_{\infty}$ can
become smaller each time we pass to a subsequence. But the special subsequence
we are about to construct has the property that the finitely many directions
composing an $\varepsilon$-net in $\mathcal{R}_{\infty}$ are stable under the
action of passing to further subsequences.

\begin{lemma}
\label{epsilon-net}Let $\varepsilon>0$ be given, and let $\left\{
\mathcal{M}_{i}^{n},g_{i}\left(  t\right)  ,O_{i},F_{i}:i\in\mathbb{N}%
\right\}  $ be a sequence that is preconverging to positive curvature. Then
there exists a finite set of directions from $\mathcal{R}_{\infty}$, say
\[
\left\{  V^{\alpha}\in\mathcal{S}_{1}^{n-1}:\alpha=1,\dots,A\right\}
\subseteq\mathcal{R}_{\infty},
\]
such that%
\[
\mathcal{R}_{\infty}\subseteq\bigcup_{\alpha=1}^{A}B_{g_{\operatorname*{can}}%
}\left(  V^{\alpha},\varepsilon\right)  .
\]
And there exists for each $a=1,\dots,A$, a particular sequence $\left\{
V_{j}^{\alpha}:j\in\mathbb{N}\right\}  $ with $\lim_{j\rightarrow\infty}%
V_{j}^{\alpha}=V^{\alpha}$ such that along a subsequence
\[
\left\{  \mathcal{M}_{i_{j}}^{n},g_{i_{j}}\left(  0\right)  ,O_{i_{j}%
},F_{i_{j}}:j\in\mathbb{N}\right\}  ,
\]
we have%
\[
\lim_{j\rightarrow\infty}\sigma_{i_{j}}\left(  V_{j}^{\alpha}\right)
=\sigma_{\infty}\left(  V^{\alpha}\right)  =\infty
\]
for all $\alpha=1,\dots,A$.
\end{lemma}

\begin{proof}
Denote by $\mathcal{R}_{\infty}^{0}$ the set of ray-like directions
corresponding to the original sequence
\[
\mathfrak{M}_{0}\doteqdot\left\{  \mathcal{M}_{i}^{n},g_{i}\left(  t\right)
,O_{i},F_{i}:i\in\mathbb{N}\right\}
\]
Then $\mathcal{R}_{\infty}^{0}\neq\emptyset$ by Lemma
\ref{R-infinity-is-nonempty}. Choose $V^{1}\in\mathcal{R}_{\infty}^{0}$, and
pass to a subsequence%
\[
\mathfrak{M}_{1}\doteqdot\left\{  \mathcal{M}_{i\left(  1,j\right)  }%
^{n},g_{i\left(  1,j\right)  }\left(  t\right)  ,O_{i\left(  1,j\right)
},F_{i\left(  1,j\right)  }:j\in\mathbb{N}\right\}
\]
along which $\sigma_{\infty}\left(  V^{1}\right)  =\infty$ is attained as the
limit%
\[
\lim_{i\left(  1,j\right)  \rightarrow\infty}\sigma_{i\left(  1,j\right)
}V_{j}^{1}=\sigma_{\infty}\left(  V^{1}\right)  =\infty
\]
for some sequence $V_{j}^{1}\rightarrow V^{1}$.

Denote by $\mathcal{R}_{\infty}^{1}\subseteq\mathcal{R}_{\infty}^{0}$ the set
of ray-like directions corresponding to the subsequence $\mathfrak{M}_{1}$.
Then $\mathcal{R}_{\infty}^{1}\neq\emptyset$ by Lemma
\ref{R-infinity-is-nonempty}. If%
\[
\mathcal{R}_{\infty}^{1}\subseteq B_{g_{\operatorname*{can}}}\left(
V^{1},\varepsilon\right)
\]
stop. Otherwise choose
\[
V^{2}\in\mathcal{R}_{\infty}^{1}\backslash B_{g_{\operatorname*{can}}}\left(
V^{1},\varepsilon\right)  ,
\]
and pass to a subsequence%
\[
\mathfrak{M}_{2}\doteqdot\left\{  \mathcal{M}_{i\left(  2,j\right)  }%
^{n},g_{i\left(  2,j\right)  }\left(  t\right)  ,O_{i\left(  2,j\right)
},F_{i\left(  2,j\right)  }:j\in\mathbb{N}\right\}
\]
along which $\sigma_{\infty}\left(  V^{2}\right)  =\infty$ is attained as the
limit%
\[
\lim_{i\left(  2,j\right)  \rightarrow\infty}\sigma_{i\left(  2,j\right)
}V_{j}^{2}=\sigma_{\infty}\left(  V^{2}\right)  =\infty
\]
for some sequence $V_{j}^{2}\rightarrow V^{2}$.

In general, denote by $\mathcal{R}_{\infty}^{\alpha}\subseteq\mathcal{R}%
_{\infty}^{\alpha-1}\subseteq\cdots\subseteq\mathcal{R}_{\infty}^{0}$ the set
of ray-like directions corresponding to the subsequence $\mathfrak{M}_{\alpha
}$. Then $\mathcal{R}_{\infty}^{\alpha}\neq\emptyset$ by Lemma
\ref{R-infinity-is-nonempty}. If%
\[
\mathcal{R}_{\infty}^{\alpha}\subseteq\bigcup_{\beta=1}^{\alpha}%
B_{g_{\operatorname*{can}}}\left(  V^{\beta},\varepsilon\right)
\]
stop. Otherwise choose
\[
V^{\alpha+1}\in\mathcal{R}_{\infty}^{\alpha}\left\backslash \bigcup_{\beta
=1}^{\alpha}B_{g_{\operatorname*{can}}}\left(  V^{\beta},\varepsilon\right)
\right.  ,
\]
and pass to a subsequence%
\[
\mathfrak{M}_{\alpha+1}\doteqdot\left\{  \mathcal{M}_{i\left(  \alpha
+1,j\right)  }^{n},g_{i\left(  \alpha+1,j\right)  }\left(  t\right)
,O_{i\left(  \alpha+1,j\right)  },F_{i\left(  \alpha+1,j\right)  }%
:j\in\mathbb{N}\right\}
\]
along which $\sigma_{\infty}\left(  V^{\alpha+1}\right)  =\infty$ is attained
as the limit%
\[
\lim_{i\left(  \alpha+1,j\right)  \rightarrow\infty}\sigma_{i\left(
\alpha+1,j\right)  }V_{j}^{\alpha+1}=\sigma_{\infty}\left(  V^{\alpha
+1}\right)  =\infty
\]
for some sequence $V_{j}^{\alpha+1}\rightarrow V^{\alpha+1}$. Since each
$\mathcal{R}_{\infty}^{\alpha}$ is contained in $\mathcal{R}_{\infty}%
^{0}\subseteq\mathcal{S}_{1}^{n-1}$ and $\mathcal{S}_{1}^{n-1}$ is compact,
this process must eventually terminate.
\end{proof}

\begin{notation}
Henceforth we shall denote the subsequence whose existence is ensured by Lemma
\ref{epsilon-net} simply by
\[
\left\{  \mathcal{M}_{i}^{n},g_{i}\left(  t\right)  ,O_{i},F_{i}%
:i\in\mathbb{N}\right\}  .
\]

\end{notation}

\bigskip

To facilitate the final step of the proof of Proposition
\ref{N_iLK-are-bounded}, we fix a length scale $\Lambda$ at which to compare
distance. To motivate our choice, consider an isosceles triangle $\Delta$ that
is symmetric about an angle $\theta\leq\theta_{0}<\pi/3$. If $\Delta$ is
embedded in Euclidean space and has side lengths $k,\ell,\ell$, then
$k\leq\ell\sqrt{2\left(  1-\cos\theta_{0}\right)  }$. In particular, given
$K\in\left(  0,\infty\right)  $, we choose $\Lambda$ depending only on $K$ to
be large enough that%
\begin{equation}
\Lambda>\frac{2K}{1-2\sqrt{2\left(  1-\cos\frac{\pi}{8}\right)  }}>0.
\label{ChooseLambda}%
\end{equation}
Combined with the simple estimate in Lemma \ref{Tangential-Radial}, this
somewhat non-intuitive choice will let us argue to a contradiction
below.\medskip

\begin{proof}
[Proof of the boundedness property]Suppose the statement is false. Then for
every $C_{j}\equiv j\in\mathbb{N}$, there exists some $L_{j}\in\left(
0,\infty\right)  $ such that for every $I\left(  j\right)  \equiv I\left(
L_{j}\right)  $, there exists some $i\left(  j\right)  \geq I\left(  j\right)
$ such that
\[
N_{i\left(  j\right)  }\left(  L_{j},K\right)  \not \subseteq B_{i\left(
j\right)  }\left(  O_{i\left(  j\right)  },j\right)  .
\]
(Recall that\ $B_{i}\left(  x,r\right)  $ denotes the open ball with center
$x\in\mathcal{M}_{i}^{n}$ and radius $r$, measured with respect to the metric
$g_{i}\left(  0\right)  $.) In particular, there exists for each $j$ some
$W_{j}\in\mathcal{S}_{1}^{n-1}$ such that%
\begin{equation}
d_{i\left(  j\right)  }\left(  O_{i\left(  j\right)  },\exp_{O_{i\left(
j\right)  }}\left(  jW_{i\left(  j\right)  }\right)  \right)  =j
\label{W-i-j-are-raylike}%
\end{equation}
and%
\begin{equation}
\exp_{O_{i\left(  j\right)  }}\left(  jW_{j}\right)  \in N_{i\left(  j\right)
}\left(  L_{j},K\right)  . \label{W-i-j-are-bad}%
\end{equation}
Notice that $W_{j}\in\mathbb{R}^{n}$ is being identified with
\[
I_{i\left(  j\right)  }\left(  W_{j}\right)  \in T_{O_{i\left(  j\right)  }%
}\mathcal{M}_{i\left(  j\right)  }^{n}.
\]
We may assume without loss of generality that%
\[
L_{j+1}\geq L_{j}+1
\]
for all $j\in\mathbb{N}$, since if $L_{\ast}\leq L^{\ast}$ then%
\[
N_{i}\left(  L_{\ast},K\right)  \subseteq N_{i}\left(  L^{\ast},K\right)
\]
for all $i\in\mathbb{N}$ and $K>0$.

We now show how to choose the $I\left(  j\right)  $. Let $\varepsilon
\in\left(  0,\pi/24\right)  $ be given. Then by Lemma \ref{epsilon-net}, there
exists a finite set of ray-like directions
\[
\left\{  V^{\alpha}:\alpha=1,\dots,A\right\}  \subseteq\mathcal{R}_{\infty}%
\]
such that%
\[
\mathcal{R}_{\infty}\subseteq\bigcup_{\alpha=1}^{A}B_{g_{\operatorname*{can}}%
}\left(  V^{\alpha},\varepsilon\right)  ;
\]
and there exist particular sequences $\left\{  V_{i}^{\alpha}\right\}  $ such
that for each $\alpha=1,\dots,A$ we have
\begin{equation}
\lim_{i\rightarrow\infty}V_{i}^{\alpha}=V^{\alpha} \label{V-alpha-i-limit}%
\end{equation}
and
\begin{equation}
\lim_{i\rightarrow\infty}\sigma_{i}\left(  V_{i}^{\alpha}\right)  =\infty.
\label{sigma-alpha-i-limit}%
\end{equation}
For each $j\in\mathbb{N}$, first choose $I^{\prime}\left(  j\right)  $ so
large that if $i\geq I^{\prime}\left(  j\right)  $, we have%
\[
\sigma_{i}\left(  V_{i}^{\alpha}\right)  \geq L_{j}%
\]
for all $a=1,\dots,A$; then choose $I\left(  j\right)  \geq I^{\prime}\left(
j\right)  $ so large that if $i\geq I\left(  j\right)  $, we have%
\[
\operatorname{sect}\left[  g_{i}\left(  x,0\right)  \right]  \geq0
\]
for all $x\in\bar{B}_{i}\left(  O_{i},3L_{j}\right)  $. This is possible by
Proposition \ref{CurvatureBoundFromBelow}.

The construction just completed yields a subsequence%
\[
\left\{  \mathcal{M}_{i\left(  j\right)  }^{n},g_{i\left(  j\right)  }\left(
0\right)  ,O_{i\left(  j\right)  },F_{i\left(  j\right)  }:j\in\mathbb{N}%
\right\}
\]
along which (\ref{W-i-j-are-raylike}) and (\ref{W-i-j-are-bad}) are satisfied
for each $i\left(  j\right)  \geq I\left(  j\right)  $. We next pass from this
to a subsequence
\[
\left\{  \mathcal{M}_{i\left(  j\left(  k\right)  \right)  }^{n},g_{i\left(
j\left(  k\right)  \right)  }\left(  0\right)  ,O_{i\left(  j\left(  k\right)
\right)  },F_{i\left(  j\left(  k\right)  \right)  }:k\in\mathbb{N}\right\}
\]
such that%
\[
\lim_{k\rightarrow\infty}W_{j\left(  k\right)  }\doteqdot W_{\infty}%
\in\mathcal{S}_{1}^{n-1}%
\]
exists. Denote by $\mathcal{R}_{\infty}^{\prime}$ the set of ray-like
directions for this subsequence. Then it may be that $\mathcal{R}_{\infty
}^{\prime}\subsetneq\mathcal{R}_{\infty}$. But by (\ref{V-alpha-i-limit}) and
(\ref{sigma-alpha-i-limit}), we still have%
\[
\left\{  V^{\alpha}:\alpha=1,\dots,A\right\}  \subseteq\mathcal{R}_{\infty
}^{\prime},
\]
and
\[
\mathcal{R}_{\infty}^{\prime}\subseteq\bigcup_{\alpha=1}^{A}%
B_{g_{\operatorname*{can}}}\left(  V^{\alpha},\varepsilon\right)  ;
\]
moreover, for each $\alpha=1,\dots,A$, we also have%
\[
\lim_{k\rightarrow\infty}V_{i\left(  j\left(  k\right)  \right)  }^{\alpha
}=V^{\alpha}%
\]
and
\[
\lim_{k\rightarrow\infty}\sigma_{i\left(  j\left(  k\right)  \right)  }\left(
V_{i\left(  j\left(  k\right)  \right)  }^{\alpha}\right)  =\infty.
\]
Observe that $W_{\infty}\in\mathcal{R}_{\infty}^{\prime}$, since we have%
\[
\sigma_{i\left(  j\left(  k\right)  \right)  }\left(  W_{j\left(  k\right)
}\right)  \geq j\left(  k\right)  \rightarrow\infty\quad\quad\text{as}%
\quad\quad k\rightarrow\infty
\]
by (\ref{W-i-j-are-raylike}). In particular, there exists some $\alpha
\in1,\dots,A$ such that%
\[
\left|  W_{\infty}-V^{\alpha}\right|  _{g_{\operatorname*{can}}}<\varepsilon.
\]

To finish the proof, choose $k$ so large that
\[
L_{j\left(  k\right)  }\geq\Lambda
\]
and that%
\[
\left|  V_{i\left(  j\left(  k\right)  \right)  }^{\alpha}-V^{\alpha}\right|
_{g_{\operatorname*{can}}}<\varepsilon
\]
and that%
\[
\sigma_{i\left(  j\left(  k\right)  \right)  }\left(  W_{j\left(  k\right)
}\right)  \geq j\left(  k\right)  \geq\Lambda
\]
and that%
\[
\left|  W_{j\left(  k\right)  }-W_{\infty}\right|  _{g_{\operatorname*{can}}%
}<\varepsilon.
\]
Then we have%
\[
\left|  W_{j\left(  k\right)  }-V_{i\left(  j\left(  k\right)  \right)
}^{\alpha}\right|  _{g_{\operatorname*{can}}}<3\varepsilon.
\]
Since
\[
\exp_{O_{i\left(  j\left(  k\right)  \right)  }}\left(  j\left(  k\right)
\cdot W_{j\left(  k\right)  }\right)  \in N_{i\left(  j\left(  k\right)
\right)  }\left(  L_{j\left(  k\right)  },K\right)
\]
and
\[
\sigma_{i\left(  j\left(  k\right)  \right)  }\left(  V_{i\left(  j\left(
k\right)  \right)  }^{\alpha}\right)  \geq L_{j\left(  k\right)  }\geq
\Lambda,
\]
we may\ apply Lemma \ref{Tangential-Radial} with $V=V_{i\left(  j\left(
k\right)  \right)  }^{\alpha}$ and $r=s=\Lambda$ to obtain the estimate%
\begin{align*}
\Lambda &  =d_{i\left(  j\left(  k\right)  \right)  }\left(  O_{i\left(
j\left(  k\right)  \right)  },\exp_{O_{i\left(  j\left(  k\right)  \right)  }%
}\left(  \Lambda W_{j\left(  k\right)  }\right)  \right) \\
&  \leq K+2\cdot d_{i\left(  j\left(  k\right)  \right)  }\left(
\exp_{O_{i\left(  j\left(  k\right)  \right)  }}\left(  \Lambda W_{j\left(
k\right)  }\right)  ,\exp_{O_{i\left(  j\left(  k\right)  \right)  }}\left(
\Lambda V_{i\left(  j\left(  k\right)  \right)  }^{\alpha}\right)  \right)  .
\end{align*}
But since%
\[
d_{i\left(  j\left(  k\right)  \right)  }\left(  \exp_{O_{i\left(  j\left(
k\right)  \right)  }}\left(  \Lambda W_{j\left(  k\right)  }\right)
,\exp_{O_{i\left(  j\left(  k\right)  \right)  }}\left(  \Lambda V_{i\left(
j\left(  k\right)  \right)  }^{\alpha}\right)  \right)  \leq2\Lambda
\leq2L_{j\left(  k\right)  },
\]
any minimizing geodesic between $\exp_{O_{i\left(  j\left(  k\right)  \right)
}}\left(  \Lambda W_{j\left(  k\right)  }\right)  $ and $\exp_{O_{i\left(
j\left(  k\right)  \right)  }}\left(  \Lambda V_{i\left(  j\left(  k\right)
\right)  }^{\alpha}\right)  $ must lie in $\bar{B}_{i\left(  j\left(
k\right)  \right)  }\left(  O_{i\left(  j\left(  k\right)  \right)
},3L_{j\left(  k\right)  }\right)  $, where the sectional curvature is
nonnegative. Hence the hinge version of the Toponogov comparison theorem
(Theorem 2.2\thinspace(B) of \cite{CE}) gives the estimate%
\begin{align*}
&  d_{i\left(  j\left(  k\right)  \right)  }\left(  \exp_{O_{i\left(  j\left(
k\right)  \right)  }}\left(  \Lambda W_{j\left(  k\right)  }\right)
,\exp_{O_{i\left(  j\left(  k\right)  \right)  }}\left(  \Lambda V_{i\left(
j\left(  k\right)  \right)  }^{\alpha}\right)  \right) \\
&  <\Lambda\sqrt{2\left(  1-\cos\left(  3\varepsilon\right)  \right)  }%
\leq\Lambda\sqrt{2\left(  1-\cos\frac{\pi}{8}\right)  }.
\end{align*}
Combining these two estimates yields $\Lambda\leq K+2\Lambda\sqrt{2\left(
1-\cos\pi/8\right)  }$, hence%
\[
0<\Lambda\leq\frac{K}{1-\sqrt{8\left(  1-\cos\frac{\pi}{8}\right)  }}<\frac
{1}{2}\Lambda
\]
by the choice we made in (\ref{ChooseLambda}). This contradiction establishes
the proposition.
\end{proof}

\section{Proof of the injectivity radius estimate\label{proof}}

The remainder of our proof of Theorem \ref{InjEstForNNCO} proceeds exactly
like the analogous part of \S 25 of \cite{H1}. But because the argument here
uses our innovations (the sets $N_{i}\left(  L,K\right)  $, for example) in an
essential way, we shall give it in detail. To prepare for this, we introduce
some notation and recall an important fact.

Let $\left(  \mathcal{M}^{n},g\right)  $ be any Riemannian manifold.

\begin{definition}
\label{DefineProperGeodesic-k-gon}If $k\in\left\{  1,2,\dots\right\}  $, a
\textbf{proper geodesic }$k$\textbf{-gon} is a collection
\[
\Gamma=\left\{  \gamma_{i}:\left[  0,\ell_{i}\right]  \rightarrow
\mathcal{M}:i=1,\dots,k\right\}
\]
of unit-speed geodesic paths between $k$ pairwise distinct vertices $p_{i}%
\in\mathcal{M}$ such that $p_{i}=\gamma_{i}\left(  0\right)  =\gamma
_{i-1}\left(  \ell_{i-1}\right)  $ for each $i$, where all indices are
interpreted modulo $k$. The \textbf{length} of a proper geodesic $k$-gon is
$L\left(  \Gamma\right)  \doteqdot\sum_{i=1}^{k}L\left(  \gamma_{i}\right)  $.
We say $\Gamma$ is a \textbf{nondegenerate proper geodesic }$k$\textbf{-gon}
if $\measuredangle_{p_{i}}\left(  -\dot{\gamma}_{i-1},\dot{\gamma}_{i}\right)
\neq0$ for each $i=1,\dots,k$; if $k=1$, we interpret this to mean $L\left(
\Gamma\right)  >0$. Finally, a \textbf{(nondegenerate) geodesic }%
$k$\textbf{-gon} is a (nondegenerate) proper geodesic $j$-gon for some
$j=1,\dots,k$.
\end{definition}

Now let $N\subset\mathcal{M}^{n}$ be a nonempty subset, and let $\Omega$
denote the space of unit-speed nondegenerate geodesic $1$-gons contained in
$N$. Let $L:\Omega\rightarrow\lbrack0,\infty)$ denote the length function, and
define $A:\Omega\rightarrow\left.  \mathcal{S}^{n-1}\mathcal{M}\right|
_{N}\times\lbrack0,\infty)$ for all unit-speed nondegenerate geodesic $1$-gons
$\alpha$ by $A\left(  \alpha\right)  \doteqdot\left(  \dot{\alpha}\left(
0\right)  ,L\left(  \alpha\right)  \right)  $. The map $A$ is injective and
induces a topology on $\Omega$ from the topology on $\left.  \mathcal{S}%
^{n-1}\mathcal{M}\right|  _{N}\times\lbrack0,\infty)$. If $N$ is compact, the
set $L^{-1}\left[  0,K\right]  \subseteq\Omega$ is compact for every
$K\in\left(  0,\infty\right)  $. If $K$ is large enough so that $L^{-1}\left[
0,K\right]  $ is nonempty, then there exists a nondegenerate geodesic $1$-gon
$\beta\in L^{-1}\left[  0,K\right]  \subset\Omega$ of minimal length. Clearly,
$\beta$ is of minimal length among all nondegenerate geodesic $1$-gons
contained in $N$; in particular, we have $L\left(  \beta\right)  =\inf
_{\alpha\in L^{-1}\left[  0,K\right]  }L\left(  \alpha\right)  =\inf
_{\alpha\in\Omega}L\left(  \alpha\right)  $.

\bigskip

\begin{proof}
[Proof of Theorem \ref{InjEstForNNCO}]Pass to a subsequence $\left\{
\mathcal{M}_{i}^{n},g_{i}\left(  t\right)  ,O_{i},F_{i}:i\in\mathbb{N}%
\right\}  $ that is preconverging to positive curvature and has the
boundedness property guaranteed by Proposition \ref{N_iLK-are-bounded}. Then
there exists $C<\infty$ such that for any $L>2$ to be chosen later, there
exists $I^{\prime}\left(  L\right)  $ such that
\[
N_{i}\left(  L,1\right)  \subseteq B_{i}\left(  O_{i},C\right)
\]
for all $i\geq I^{\prime}\left(  L\right)  $. By Proposition
\ref{CurvatureBoundFromBelow}, there exist $I\left(  L\right)  \geq I^{\prime
}\left(  L\right)  $ and $\eta>0$ such that for all $i\geq I\left(  L\right)
$, we have%
\begin{equation}
\inf\left\{  \operatorname{sect}\left(  g_{i}\left(  x,0\right)  \right)
:x\in\bar{B}_{i}\left(  O_{i},C+2\right)  \right\}  \geq\eta.
\label{BumpAtOrigin}%
\end{equation}

Suppose the theorem is false. Then there exists $i_{0}$ such that
$\operatorname*{inj}{}_{g_{i}\left(  0\right)  }\left(  O_{i}\right)  <1$ for
all $i\geq i_{0}$. So there exists for each $i\geq i_{0}$ a nondegenerate
geodesic $2$-gon $\alpha_{i}$ based at $O_{i}$ and of length $<2$, hence
contained in $B_{i}\left(  O_{i},1\right)  $. By Lemma
\ref{N_iLK-contains-K-ball}, $\bar{B}_{i}\left(  O_{i},1\right)  \subseteq
N_{i}\left(  L,1\right)  $. So by a standard shortening argument, there exists
for each $i\geq i_{0}$ a nondegenerate geodesic $1$-gon $\tilde{\alpha}_{i}$
based at $O_{i}$, contained in $N_{i}\left(  L,1\right)  $, and such that
$\operatorname{length}_{g_{i}\left(  0\right)  }\tilde{\alpha}_{i}%
<\operatorname{length}_{g_{i}\left(  0\right)  }\alpha_{i}<2$. By Lemma
\ref{N_iLK-are-compact}, $N_{i}\left(  L,1\right)  $ is compact. So there
exists for each $i\geq i_{0}$ a shortest element $\beta_{i}$ in the set of all
nondegenerate geodesic $1$-gons contained in $N_{i}\left(  L,1\right)  $. Each
$\beta_{i}$ is smooth except perhaps at its base $\beta_{i}\left(  0\right)
=\beta_{i}\left(  \ell_{i}\right)  $, where
\[
\ell_{i}\doteqdot\operatorname{length}_{g_{i}\left(  0\right)  }\beta_{i}%
\leq\operatorname{length}_{g_{i}\left(  0\right)  }\tilde{\alpha}_{i}<2.
\]

We first consider the (easier) case that there exists a subsequence for which
$\beta_{i}$ is smooth at $\beta_{i}\left(  0\right)  =\beta_{i}\left(
\ell_{i}\right)  $. By Lemma \ref{R-infinity-is-nonempty}, the set
$\mathcal{R}_{\infty}$ for this subsequence is nonempty. Hence definition
(\ref{sigma-infinity}) implies that for every $L$ and every $J\in\mathbb{N}$,
there exists some $i\left(  I\left(  L\right)  ,J\right)  \geq\max\left\{
I\left(  L\right)  ,J\right\}  $ and some $V_{i}\in\mathcal{S}_{1}^{n-1}$ such
that $\sigma_{i}\left(  V_{i}\right)  \geq L$. Let%
\[
y_{i}\doteqdot\exp_{O_{i}}\left(  LV_{i}\right)  ,
\]
and define%
\[
S_{i}\doteqdot d_{i}\left(  y_{i},\beta_{i}\right)  .
\]
Since $\beta_{i}\subset N_{i}\left(  L,1\right)  $ is compact, there are
$W_{i}\in\mathcal{S}_{1}^{n-1}$ and $t_{i}\in\left[  0,\sigma_{i}\left(
W_{i}\right)  \right]  \backslash\left\{  \infty\right\}  $ such that%
\[
x_{i}\doteqdot\exp_{O_{i}}\left(  t_{i}W_{i}\right)  \in\beta_{i}%
\]
is the point on $\beta_{i}$ closest to $y_{i}$, so that $d_{i}\left(
x_{i},y_{y}\right)  =S_{i}$. Since $x_{i}\in N_{i}\left(  L,1\right)  $, the
definition of $N_{i}\left(  L,1\right)  $ implies that%
\[
S_{i}=d_{i}\left(  x_{i},y_{i}\right)  \geq L-1>1.
\]
Let $\gamma_{i}$ be a minimal unit-speed geodesic from $x_{i}$ to $y_{i}$.
Note in particular that $\operatorname{length}_{g_{i}\left(  0\right)  }%
\gamma_{i}=S_{i}>1$ and that%
\begin{equation}
\left.  \gamma_{i}\right|  _{\left[  0,1\right]  }\subset B_{i}\left(
O_{i},C+1\right)  . \label{gamma-in-bump}%
\end{equation}
Since $\beta_{i}$ is smooth, we can apply the first variation formula to
conclude that $\dot{\beta}_{i}\perp\dot{\gamma}_{i}$ at $x_{i}$, where
$\dot{\beta}_{i},\dot{\gamma}_{i}$ denote the unit tangent vectors of
$\beta_{i},\gamma_{i}$ respectively. Let $X_{i}$ be the unit vector field that
results from parallel translation of $\dot{\beta}_{i}$ along $\gamma_{i}$ from
$x_{i}$, and define the cutoff function%
\[
f_{i}\left(  s\right)  =\left\{
\begin{array}
[c]{cl}%
1 & \text{if\quad}0\leq s\leq1\\
\left(  S_{i}-s\right)  /\left(  S_{i}-1\right)  & \text{if\quad}1<s\leq S_{i}%
\end{array}
.\right.
\]
Then the minimality of $\gamma_{i}$ implies that the second-variation index
form $\mathcal{I}$ in the direction $f_{i}X_{i}$ is nonnegative:%
\[
0\leq\mathcal{I}\equiv\mathcal{I}\left(  f_{i}X_{i},f_{i}X_{i}\right)
\doteqdot\int_{\gamma_{i}}\left(  \left|  \nabla_{\dot{\gamma}_{i}}\left(
f_{i}X_{i}\right)  \right|  ^{2}-\left\langle R\left(  \dot{\gamma}_{i}%
,f_{i}X_{i}\right)  \left(  f_{i}X_{i}\right)  ,\dot{\gamma}_{i}\right\rangle
\right)  \,ds.
\]
But (\ref{BumpAtOrigin}) and (\ref{gamma-in-bump}) imply that all sectional
curvatures are bounded below by $\eta>0$ along $\left.  \gamma_{i}\right|
_{\left[  0,1\right]  }$. And assumption (\ref{AlmostNonNegative}) implies
that all sectional curvatures are bounded below by $-\delta_{i}\nearrow0$
throughout $\mathcal{M}_{i}^{n}$. Hence we can estimate%
\begin{align*}
\mathcal{I}  &  =-\int_{0}^{1}\left\langle R\left(  \dot{\gamma}_{i}%
,X_{i}\right)  X_{i},\dot{\gamma}_{i}\right\rangle \,ds+\int_{1}^{S_{i}%
}\left(  \left(  df_{i}\left(  \dot{\gamma}_{i}\right)  \right)  ^{2}%
-f_{i}^{2}\left\langle R\left(  \dot{\gamma}_{i},X_{i}\right)  X_{i}%
,\dot{\gamma}_{i}\right\rangle \right)  \,ds\\
&  \leq-\eta+\frac{1}{S_{i}-1}+\delta_{i}\frac{S_{i}-1}{3}.
\end{align*}
Now choose $L$ so large that $S_{i}\geq L-1$ satisfies $1/\left(
S_{i}-1\right)  \leq\eta/3$. Then choose $J$ so large that for all $i\geq J$,
we have $\delta_{i}\left(  S-1\right)  \leq\eta$. Thus for $i=i\left(
I\left(  L\right)  ,J\right)  $, we get $\mathcal{I}\leq-\eta/3<0$. This
contradicts the minimality of $\gamma_{i}$ and proves the theorem in this case.

Now we consider the case that there exists $i_{1}\geq i_{0}$ such that
$\beta_{i}$ fails to be smooth at $\beta_{i}\left(  0\right)  =\beta
_{i}\left(  \ell_{i}\right)  $ for all $i\geq i_{1}$. It is a standard fact
that\ for any complete Riemannian manifold $\left(  \mathcal{M}^{n},g\right)
$ with sectional curvatures bounded above by $\kappa>0$, any points
$p,q\in\mathcal{M}^{n}$, any geodesic path $\gamma$ from $p$ to $q$ of length
less than $\pi/\sqrt{\kappa}$, and any points $\tilde{p},\tilde{q}$
sufficiently near $p,q$ respectively, there exists a unique geodesic
$\tilde{\gamma}$ from $\tilde{p}$ to $\tilde{q}$ that is close to $\gamma$.
Consider a variation moving $\beta_{i}\left(  \ell_{i}\right)  $ in the
direction $\dot{\beta}_{i}\left(  0\right)  $ and observe that the first
variation in this direction is strictly negative: $\left\langle \dot{\beta
}_{i}\left(  \ell_{i}\right)  ,\dot{\beta}_{i}\left(  0\right)  \right\rangle
-\left\langle \dot{\beta}_{i}\left(  0\right)  ,\dot{\beta}_{i}\left(
0\right)  \right\rangle <0$. Now we have $\operatorname{sect}\left(
g_{i}\left(  0\right)  \right)  \leq1$ by assumption (\ref{AlmostNonNegative}%
), and $\operatorname{length}_{g_{i}\left(  0\right)  }\beta_{i}=\ell_{i}<1$
by hypothesis. It follows that there exists a nondegenerate geodesic $1$-gon
$\tilde{\beta}_{i}$ with%
\[
\operatorname{length}_{g_{i}\left(  0\right)  }\tilde{\beta}_{i}%
<\operatorname{length}_{g_{i}\left(  0\right)  }\beta_{i},
\]
and such that
\[
\tilde{\beta}_{i}\left(  0\right)  \in\beta_{i}\subset N_{i}\left(
L,1\right)
\]
and%
\[
\tilde{\beta}_{i}\subset B_{i}\left(  O_{i},C+1\right)  .
\]
$\tilde{\beta}_{i}$ may not be smooth at its base $\tilde{\beta}_{i}\left(
0\right)  $ either, but it is smooth everywhere else.

By our choice of $\beta_{i}$, it must be that $\tilde{\beta}_{i}$ does not lie
entirely in $N_{i}\left(  L,1\right)  $. Hence there must exist a point
$z_{i}$ on $\tilde{\beta}_{i}$ but not in $N_{i}\left(  L,1\right)  $. Choose
$W_{i}\in\mathcal{S}_{1}^{n-1}$ and $t_{i}\in\left[  0,\sigma_{i}\left(
W_{i}\right)  \right]  \backslash\left\{  \infty\right\}  $ such that%
\[
z_{i}\doteqdot\exp_{O_{i}}\left(  t_{i}W_{i}\right)  \in\tilde{\beta}%
_{i}\backslash N_{i}\left(  L,1\right)  .
\]
By definition of $N_{i}\left(  L,1\right)  $, there exist some $V_{i}%
\in\mathcal{S}_{1}^{n-1}$ with $\sigma_{i}\left(  V_{i}\right)  \geq L$, some
$r_{i}\in\left[  0,\sigma_{i}\left(  V_{i}\right)  \right]  \backslash\left\{
\infty\right\}  $, and some $s_{i}\in\left[  0,t_{i}\right]  $ such that%
\begin{equation}
d_{i}\left(  \exp_{O_{i}}\left(  r_{i}V_{i}\right)  ,\exp_{O_{i}}\left(
s_{i}W_{i}\right)  \right)  <r_{i}-1. \label{NotInN_iLK}%
\end{equation}
Define%
\[
y_{i}\doteqdot\exp_{O_{i}}\left(  LV_{i}\right)  ,
\]
and let $\zeta_{i}$ denote the geodesic%
\[
\zeta_{i}:\left[  0,t_{i}\right]  \rightarrow\mathcal{M}_{i}^{n},\quad
\quad\quad\quad\zeta_{i}:\tau\mapsto\exp_{O_{i}}\left(  \tau W_{i}\right)  .
\]
We claim that $z_{i}=\zeta_{i}\left(  t_{i}\right)  $ is the point on
$\zeta_{i}$ closest to $y_{i}$. To see this, first note that the closest point
is not $O_{i}=\zeta_{i}\left(  0\right)  $, since (\ref{NotInN_iLK}) implies
that%
\begin{align}
d_{i}\left(  y_{i},\zeta_{i}\left(  s_{i}\right)  \right)   &  \leq
d_{i}\left(  y_{i},\exp_{O_{i}}\left(  r_{i}V_{i}\right)  \right)
+d_{i}\left(  \exp_{O_{i}}\left(  r_{i}V_{i}\right)  ,\zeta_{i}\left(
s_{i}\right)  \right) \nonumber\\
&  <\left(  L-r_{i}\right)  +\left(  r_{i}-1\right)  =d_{i}\left(  y_{i}%
,O_{i}\right)  -1. \label{y-is-close-to-zeta}%
\end{align}
If the closest point is an interior point, say $\zeta_{i}\left(  \tau
_{i}\right)  $ for some $\tau_{i}\in\left(  0,t_{i}\right)  $, let $\xi_{i}$
be a minimal geodesic from $\zeta_{i}\left(  \tau_{i}\right)  $ to $y_{i}$.
Because $z_{i}\in\tilde{\beta}_{i}\subset B_{i}\left(  O_{i},C+1\right)  $,
and all sectional curvatures are bounded below by $\eta>0$ in $B_{i}\left(
O_{i},C+2\right)  $, a second variation argument (like the one above) along
$\xi_{i}$ will yield a contraction. This proves that the closest point to
$y_{i}$ on $\zeta_{i}$ cannot be an interior point. Hence the only possibility
is that the closest point to $y_{i}$ along $\zeta_{i}$ is its other endpoint
$z_{i}=\zeta_{i}\left(  t_{i}\right)  $. This proves the claim. (Note that
applying this argument along segments proves the stronger fact that the
function $\tau\mapsto d_{i}\left(  y_{i},\zeta_{i}\left(  \tau\right)
\right)  $ is monotone decreasing for $\tau\in\left[  0,t_{i}\right]  $.)

By the claim and (\ref{y-is-close-to-zeta}), we have%
\[
d_{i}\left(  y_{i},z_{i}\right)  \leq d_{i}\left(  y_{i},\zeta_{i}\left(
s_{i}\right)  \right)  <L-1.
\]
But since $\tilde{\beta}_{i}\left(  0\right)  \in N_{i}\left(  L,1\right)  $,
we have%
\[
d_{i}\left(  y_{i},\tilde{\beta}_{i}\left(  0\right)  \right)  \geq L-1.
\]
Hence the closest point to $y_{i}$ on $\tilde{\beta}_{i}$ is not its base
$\tilde{\beta}_{i}\left(  0\right)  $. In particular, $\tilde{\beta}_{i}$ is
smooth at its closest point to $y_{i}$. So we can construct a
length-minimizing geodesic $\tilde{\gamma}_{i}$ from $\tilde{\beta}_{i}$ to
$y_{i}$ and apply a second variation argument along it, exactly as in the
first case. Because $\left.  \tilde{\gamma}_{i}\right|  _{\left[  0,1\right]
}\subset B_{i}\left(  O_{i},C+2\right)  $, where the sectional curvatures are
bounded from below by $\eta>0$, this argument leads to a contradiction, just
as before. This finishes the proof.\medskip
\end{proof}

\end{document}